\documentclass[sn-mathphys-num]{sn-jnl}

\usepackage{graphicx}%
\usepackage{multirow}%
\usepackage{amsmath,amssymb,amsfonts}%
\usepackage{amsthm}%
\usepackage{mathrsfs}%
\usepackage[title]{appendix}%
\usepackage{xcolor}%
\usepackage{textcomp}%
\usepackage{manyfoot}%
\usepackage{booktabs}%
\usepackage{algorithm}%
\usepackage{algorithmicx}%
\usepackage{algpseudocode}%
\usepackage{listings}%

\newtheorem{conjecture}{Conjecture}

\theoremstyle{thmstyleone}%
\theoremstyle{thmstyletwo}%
\newtheorem{example}{Example}%

\theoremstyle{thmstylethree}%
\newtheorem{definition}{Definition}%

\usepackage[acronym]{glossaries}
\glsdisablehyper
\newacronym{dfo}{DFO}{derivative-free optimization}
\usepackage{booktabs}

\usepackage{color}

\begin{document}

\title[On the Relationship between \(\Lambda\)-poisedness in Derivative-Free Optimization and Outliers in Local Outlier Factor]{On the Relationship between \(\Lambda\)-poisedness in Derivative-Free Optimization and Outliers in Local Outlier Factor}


\author[1]{\fnm{Qi} \sur{Zhang}}\email{zqaptx@163.com}

\author[2]{\fnm{Pengcheng} \sur{Xie}}\email{pxie@lbl.gov}
\equalcont{These authors contributed equally to this work.}


\affil*[1]{\orgdiv{School of Electrical Engineering}, \orgname{Xi'an Jiaotong University}, \orgaddress{\street{Taiyilu Avenue}, \city{Xi'an}, \postcode{710049}, \state{Shaanxi Province}, \country{China}}}

\affil[2]{\orgdiv{Applied Mathematics and Computational Research Division}, \orgname{Lawrence Berkeley National Laboratory}, \orgaddress{\street{1 Cyclotron Road}, \city{CA}, \postcode{94720}, \state{State}, \country{United States}}}



\abstract{Derivative-free optimization (DFO) is a method that does not require the calculation of gradients or higher-order derivatives of the objective function, making it suitable for cases where the objective function is non-differentiable or the computation of derivatives is expensive. This communication discusses the importance of \(\Lambda\)-poisedness in DFO and the outliers detected by the Local Outlier Factor (LOF) on the optimization process. We discuss the relationship between \(\Lambda\)-poisedness in derivative-free optimization and outliers in local outlier factor.}

\keywords{Derivative-free optimization, well-poisedness}


\pacs[MSC Classification]{90C56, 41A10, 65K05, 90C300}

\maketitle

\section{Introduction}\label{sec1}

Most optimization methods for unconstrained optimization require the use of the derivative of the objective function. However, in practical applications, the objective function is often expensive to evaluate, and its derivatives are not available. Derivative-free optimization methods address the optimization problem
\begin{equation}
    \label{eq:nlp}
    \min_{x \in  \Re^n}  f(x)
\end{equation}
by relying solely on function evaluations without the need for derivative information.

Examples of such problems are found in various fields, such as reinforcement learning~\cite{Qian_Yu_2021}, hyperparameter tuning~\cite{Ghanbari_Scheinberg_2017}, particle physics~\cite{Eldred_Etal_2023}, etc. One can see the monographs~\cite{Conn_Scheinberg_Vicente_2009,Audet_Hare_2017}, the survey papers~\cite{Powell_1998,Custodio_Scheinberg_Vicente_2017,Larson_Menickelly_Wild_2019}, and our previous work \cite{xie2023dfoto,xie2023least,ragonneau2024optimal,xie2023linesearch} for more details.

This communication is organized as follows. 
Section~\ref{sec:poisedness-lof} then gives the concept of LOF and the bad point under  $\Lambda$-poisedness. Section~\ref{sec:numerical-examples} shows the example we used and briefly gives some numerical results. Section~\ref{sec:poisedness-versus-lof} carries out a series of numerical experiments. Section~\ref{sec:conclusion} concludes this communication.

\section{LOF and $\Lambda$-poisedness}
\label{sec:poisedness-lof}

On the one hand, $\Lambda$-poisedness is an essential criterion of well-poisedness of an interpolation set. Certainly, it can also be used for identifying those bad interpolation points, while the bad point under $\Lambda$-poisedness is sometimes not that obvious or hard to explain directly. On the other hand, as a famous outlier identification criterion, LOF can be used to select those points far away from others in the same point set. Compared with $\Lambda$-poisedness, the geometrical meaning of LOF is single, clear, and more explicit. Therefore, in this section, we try to rediscover the meaning of $\Lambda$-poisedness from the aspect of LOF. With the help of LOF, it's hopeful to gain a deeper understanding of $\Lambda$-poisedness from their differences and similarities.

\subsection{LOF}

LOF is a density-based outliers detection algorithm. It was first proposed in 2000 by Breunig et al \cite{breunig_lof_2000} to overcome the non-local problem of previous detection algorithms. This algorithm can give the outlierness of a data point quantitatively. The larger the LOF of a data point, the more likely it is an outlier. LOF algorithm has been used in various fields successfully, such as hyperspectral imagery noisy label detection \cite{tu_hyperspectral_2018}, electricity theft detection \cite{peng_electricity_2021} and even in switching median filter \cite{wang_efficient_2011}.

In this subsection we will give the definition of LOF and show the explicit expression of the LOF of a data point under certain special conditions. First of all, we begin with the definition of distance used in LOF.

\begin{definition}[$k$-distance of an object $p$]
For any positive integer $k$, the $k$-distance of object $p$, denoted as $k$-$\operatorname{distance}(p)$, is defined as the distance $d(p, o)$ between $p$ and an object $o \in D$ such that:

(i) for at least $k$ objects $o^{\prime} \in D \backslash\{p\}$ it holds that $d\left(p, o^{\prime}\right) \leq d(p, o)$, and

(ii) for at most $k-1$ objects o' $\in D \backslash\{p\}$ it holds that $d\left(p, o^{\prime}\right)<d(p, o)$.
\end{definition}

$k$-distance of a data point, called center point afterwards, provides a clear boundary of certain scope of a data point. Points in this scope are treated as neighbors of that center point and the scope is called its neighborhood. Here shows the strict definition of the scope above.

\begin{definition} [$k$-distance neighborhood of an object $p$]
Given the $k$-distance of $p$, the $k$-distance neighborhood of $p$ contains every object whose distance from $p$ is not greater than the $k$-distance, i.e. $N_{k\operatorname{-distance}(p)}(p)=\{q \in D \backslash\{p\} \mid d(p, q) \leq k$-$\operatorname{distance}(p)\}$, noted as $N_k(p)$.
These objects $q$ are called the $k$-nearest neighbors of $p$.
\end{definition}

Obviously, the $k$-distance of a center point can directly used for measuring the local density of this center point. However, it still lacks of some global information of the whole poised set and is not even a distance. Therefore, we still need the following definition.

\begin{definition}[reachability distance of an object $p$ w.r.t. object $o$]
Let $k$ be a natural number. The reachability distance of object $p$ with respect to object $o$ is defined as
\[
\text{\rm reach-dist}_k(p, o)=\max \{k\operatorname{-distance}(o), d(p, o)\}.
\]

\end{definition}

The reachability distance of two data points depends not only on the absolute distance of these two points but also on the neighborhood of the second point. All the points in $k$-distance neighborhood of point $o$ will be pushed to the edge of its $k$-distance neighborhood. This indicates that the distance between points $p$ and $o$ may differ from the distance between points $o$ and $p$. 

\begin{figure}[htbp]
    \centering
    \includegraphics[scale=0.7]{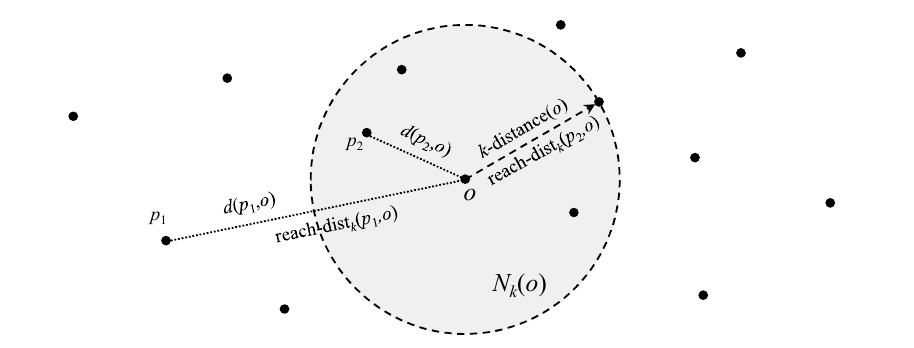}
    \caption{$k$-distance and reachability distance, with $k$=4. The $\text{\rm reach-dist}_k$ of two points is controlled by the $k$-distance of the second point, instead of simply the distance between these two points.}
\end{figure} 

With the definitions above, we can finally define the density and outlier factor of a point in LOF.

\begin{definition}[local reachability density of an object $p$] 
The local reachability density of $p$ is defined as
$$
\text{lrd}_{\text{MinPts}}(p)=1 {\Large\slash} \left(\frac{\sum_{o \in N_{\text{MinPts }}(p)} \text{ reach-dist }_{\text{MinPts}}(p, o)}{\left|N_{\text{MinPts}}(p)\right|}\right)
$$
where $\text{MinPts}$ is the parameter specifying a minimum number of objects, like some typical density-based clustering algorithm.
\end{definition}

\begin{definition}[LOF]
The (local) outlier factor of $p$ is defined as
\begin{equation}
\label{eq:lof}
\text{LOF}_{\text{MinPts}}(p)=\frac{1}{{\left|N_{\text{MinPts}}(p)\right|}}\sum_{o \in N_{\text{MinPts}}(p)} \frac{\text{lrd}_{\text{MinPts}}(o)}{\text{lrd}_{\text{MinPts}}(p)}
\end{equation}
\end{definition}

When used in an interpolation problem, the expression of LOF could be reformed in a more clear and direct way. Supposing that we have a poised set $Y=\{y^0,y^1,\cdots,y^p\}$, taking $MinPts=p-1$, for some point $y^\tau,y^j\in Y$, we will have:
$$
\textrm{reach-dist}_k(y^\tau,y^j)=\underset{y^\nu\in Y}{\arg\max}\ d(y^j,y^\nu)
$$
Then for local reachable distance of point $y^\tau$:
$$
\text{lrd}_k(y^\tau) = \frac{p-1}{\displaystyle\sum_{j=1,j\neq\tau}^p\underset{y^\nu\in Y}{\arg\max}\ d(y^j,y^\nu)}
$$
Finally, we can obtian the LOF of $y^\tau$ as
$$
\begin{aligned}
{\rm LOF}_{k}&=\frac{1}{p-1}\sum_{i=1,i\neq \tau}^{p}\frac{\displaystyle\sum_{j=1,j\neq\tau}^p\underset{y^\nu\in Y}{\arg\max}\ d(y^j,y^\nu)}{\displaystyle\sum_{j=1,j\neq i}^p\underset{y^\nu\in Y}{\arg\max}\ d(y^j,y^\nu)}. 
\end{aligned}
$$

Besides, the expression of LOF with $MinPts=1$ also can be rewritten in a simple way. However for other conditions, it's difficult to obtian an explicit expression of LOF directly like this.

\subsection{Bad point under $\Lambda$-poisedness}

To rediscover the property of $\Lambda$-poisedness through LOF directly, one effective way is comparing the bad points selected by LOF with those under the sence of $\Lambda$-poisedness. We use the following criterion to choose bad point with respect to $\Lambda$-poisedness.

\begin{sloppypar}
\begin{definition}[{Bad point under $\Lambda$-poisedness}]

Let $\Lambda>0$ and a set $B \in \Re^n$ be given. Let $\phi=\left\{\phi_0(x), \phi_1(x), \ldots, \phi_p(x)\right\}$ be a basis in $\mathcal{P}_n^d$. A point $y^k$ in poised set $Y=\left\{y^0, y^1, \ldots, y^p\right\}$ is said to be bad under  $\Lambda$-poised in $B$ (in the interpolation sense) if and only if for the basis of Lagrange polynomials associated with $Y$
$$
\max _{x \in B}\left|\ell_k(x)\right| \geq \max _{x \in B}\left|\ell_i(x)\right|
$$
for $\forall i, 0 \leq i \leq p, i \neq k$.

\end{definition}
\end{sloppypar}

By testing the maximum value of each interpolation point's Lagrange polynomial in certain region, we can easily get the point with the largest maximum value. This is the bad point under the sence of $\Lambda$-poisedness. Besides, to make the comparison between LOF and $\Lambda$-poisedness as reasonable as possible, we use the ball determined by the $k$-distance of an interpolation point, noted as $\mathcal{Y}_k^i$, as the region above.

\section{Numerical Examples}
\label{sec:numerical-examples}

In this section, we will give a detailed introduction to the 2-dimensional example we mainly focus on in the following parts. There are mainly three reasons why this example is selected. The first reason is that it must be representative. Different from 1-dimensional, 2-dimensional can give us a partial glance at the condition in higher dimensions. What's more, results in 2-dimensional are much easier to visualize. An explicit image can always provide us with a clearer view of the central question. Last but not least, it's also convenient for demonstrating.

Here we give the basic example of an interpolation in 2-dimensional case.

\begin{example}[{A 2-dimensional example}]
\label{example-1}
Consider five points $y^i=(x_1^i,x_2^i)\in\mathcal{Y} \subseteq\Re^2$, $i=1,2,\cdots,N$ and a function $f(y):\Re^2\to\Re$, where $\mathcal{Y}\subseteq\mathcal{D}(f)$ the domain of $f$. At this time, all five points are the same, without priority. Without losing generality, two of them are fixed at $\Tilde{y}^1$ and $\Tilde{y}^2$, noted as $y^1$ and $y^2$. The other three still are free. Then another two points are randomly selected from these three, noted as $y^3$ and $y^4$. These two points are treated as movable points, and they can be any point in $\mathcal{Y}$, except $\Tilde{y}^1$ and $\Tilde{y}^2$. However, after their position is located, noted as $\Tilde{y}^3$ and $\Tilde{y}^4$, then they become non-movable. That is to say, as seen from the last point $y^5$, all the positions of other points are fixed. We mainly focus on $y^5$. The bad points under  $\Lambda$-poisedness and LOF, with $y^5$ moving through $\mathcal{Y}\backslash\{\Tilde{y}^1,\Tilde{y}^2,\Tilde{y}^3,\Tilde{y}^4\}$, can be selected by
$$
\begin{aligned}
i = & \mathop{\arg\max}\limits_{0\le i\le 5}{\max_{y\in\mathcal{Y}_k^i}{\left\lvert\ell_i(y)\right\rvert}},\\
j = & \mathop{\arg\max}\limits_{0\le i\le 5}{\mathrm{LOF}_k(y^i)}. 
\end{aligned}
$$

\end{example}

Based on the example above, we will give a basic numerical result for demonstration. The Lagrange interpolation formula of $y^i$ can be written as
$$
\ell_{i}(y)=y^\top A_i y+b_iy+c_i,i=1,2,\cdots,5.
$$
By solving the following equation, we can obtain the coefficients of all formulas, 
$$
\sum_{k}f(y_i)\ell_k(y_i)=f(y_i),i=1,2,\cdots,5.
$$
 To illustrate Lagrange polynomials in $\Re^2$, consider interpolating the cubic function given by
$$
f\left(x_1, x_2\right)=x_1+x_2+2 x_1^2+3 x_2^3
$$
at the six interpolating points $y^0=(0,0), y^1=(1,0), y^2=(0,1), y^3=(2,0), y^4=(1,1)$, and $y^5=(0,2)$. Clearly, $f\left(y^0\right)=0, f\left(y^1\right)=3, f\left(y^2\right)=4, f\left(y^3\right)=10, f\left(y^4\right)=7$, and $f\left(y^5\right)=26$. It is easy to see that the corresponding Lagrange polynomials $\ell_j\left(x_1, x_2\right), j=$ $0, \ldots, 5$, are given by
$$
\left\{
\begin{aligned}
& \ell_0\left(x_1, x_2\right)=1-\frac{3}{2} x_1-\frac{3}{2} x_2+\frac{1}{2} x_1^2+\frac{1}{2} x_2^2+x_1 x_2, \\
& \ell_1\left(x_1, x_2\right)=2 x_1-x_1^2-x_1 x_2, \\
& \ell_2\left(x_1, x_2\right)=2 x_2-x_2^2-x_1 x_2, \\
& \ell_3\left(x_1, x_2\right)=-\frac{1}{2} x_1+\frac{1}{2} x_1^2, \\
& \ell_4\left(x_1, x_2\right)=x_1 x_2, \\
& \ell_5\left(x_1, x_2\right)=-\frac{1}{2} x_2+\frac{1}{2} x_2^2 .
\end{aligned}
\right. 
$$

With all the Lagrange interpolation formulas known, we can get their max value in $\mathcal{Y}$ by solving a constrained optimization problem, However, based on the condition above, the actual value we want is the maximum of $l_i(y),y\in\mathcal{Y}_k^i$, and the $k$-distances of each interpolation point with respect to others are still unknown. To get the $k$-distances of these points, an essential parameter of LOF, $k$, is needed. Considering the fact that only one bad point will be selected by $\Lambda$-poisedness at once, we choose $k$ to be $5-2=3$. This will make LOF more likely to select the point farthest away from all the others as an outlier.

\begin{figure}[htbp]
    \centering
    \begin{minipage}[t]{0.45\linewidth}
        \centering
        {\includegraphics[width=0.8\textwidth]{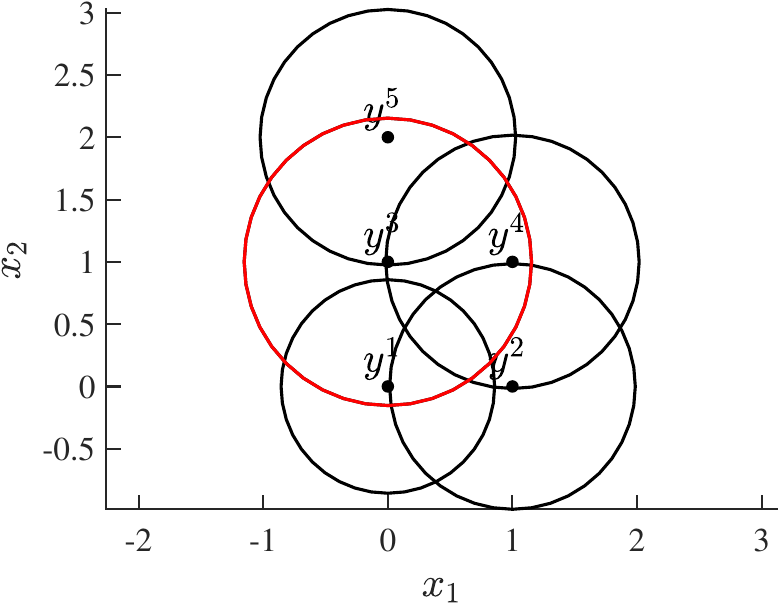}}
      
    \end{minipage}\ 
    \begin{minipage}[t]{0.45\linewidth}
        \centering
    {\includegraphics[width=0.8\textwidth]{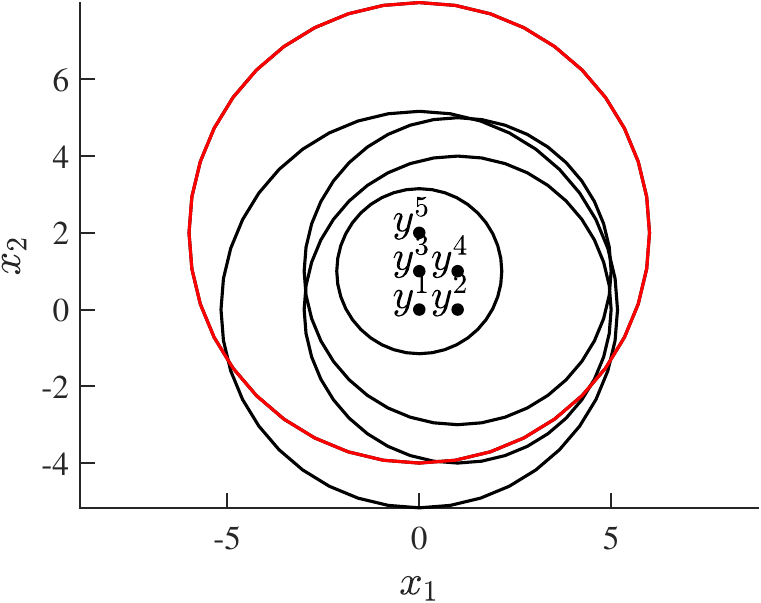}}
    \end{minipage}
      \caption{LOF (left-hand side) versus  {$\Lambda$-poisedness} (right-hand side) of Example \ref{example-1}. The size of the circle centered at each point reflects the value of the corresponding factor of that point, and the largest is highlighted in red.}
        \label{fig:numeric-lof-lambda}
\end{figure}

Table \ref{tab:numeric-lof-lambda} shows the value of $k$-distance, $\Lambda$-poisedness, and the LOF, which can be calculated through (\ref{eq:lof}) directly. From the values in the table, it's easy to tell that the bad point is $y^3$ under LOF, while $y^5$ under $\Lambda$-poisedness. What's more, though we call $y^3$ the "bad point" under LOF, it's not a real or typical outlier. Strictly, none of these five points is an outlier, because all their LOF value is near 1. For convenience and numerical analysis, we call the point with the largest LOF value and larger than 1.2 as the bad point under LOF. Fig.  \ref{fig:numeric-lof-lambda} further shows the spatial relation of these five points. The LOF of each point is represented by the size of the circle centered at the point in the left panel of Fig.  \ref{fig:numeric-lof-lambda}, and the same for $\Lambda$-poisedness in the right panel of Fig.  \ref{fig:numeric-lof-lambda}. The points with the largest LOF or $\Lambda$-poisedness value are highlighted in red. 

\begin{table}[htbp]
	\centering
 \caption{The $k$-distance, LOF and $\lambda$-poisedness value of $y^1,y^2,y^3,y^4,y^5$ in Example \ref{example-1}}
	\label{tab:numeric-lof-lambda}
 \setlength{\tabcolsep}{15pt} 
	\begin{tabular}{llllll}
    \toprule
		                   & $y^1$ & $y^2$ & $y^3$ & $y^4$ & $y^5$ \\
    \midrule
		$k$-distance          & 1.4142 & 1.4142 & 1.0000 & 1.4142 & 2.0000\\
	  LOF                   & 0.8570 & 0.9858 & 1.1532 & 1.0169 & 1.0257\\
		$\Lambda$-poisedness  & 5.1666 & 4.0000 & 2.1505 & 4.0000 & 6.0000\\
    \bottomrule
	\end{tabular}
	
\end{table}

\section{$\Lambda$-poisedness versus LOF}
\label{sec:poisedness-versus-lof}

Based on the Example \ref{example-1}, in this subsection we will carry out a series of numerical experiments. By observing and analyzing the experimental results, we find some unexpected phenomena. After summarizing and generalizing these phenomena, some meaningful conclusions are drawn. Combined with the theory about $\Lambda$-poisedness and LOF, we try to analyze these conclusions theoretically.

In order to carry out the numerical experiments, we first need to define the observation region and precision of the experiments. Considering the computational burden and observation requirements, we limit the observation area to a two-dimensional region $\Tilde{\mathcal{Y}}$ of $[-5,5]\times[-5,5]$, and take the points in the region in steps of 1 to form the set $\mathcal{Y}$. 

Then we further choose the positions of $y^1$ and $y^2$ as $(-1,0)$ and $(1,0)$. There is no special requirement for the location of their position, only for the convenience of observation, so the center of the region is chosen. For the remaining three points, as mentioned before, they can move freely within the region $\mathcal{Y}$. With the changing of the position of $y^3$ and $y^4$, We mainly focus on the state of the last point $y^5$. By traversing all the optional positions within the region, the state of y5 will be recorded. It may be an outlier under LOF or a bad point under $\Lambda$-poisedness, or, of course, both or neither.

Among all these states, we are mainly interested in the bad points under $\Lambda$-poisedness, and the position which is both a bad point under $\Lambda$-poisedness and an outlier under LOF.

\subsection{{Outlier Trap}}

\begin{figure}[htbp]
    \centering
    \begin{minipage}[t]{0.49\linewidth}
        \centering
        {\includegraphics[width=0.8\textwidth]{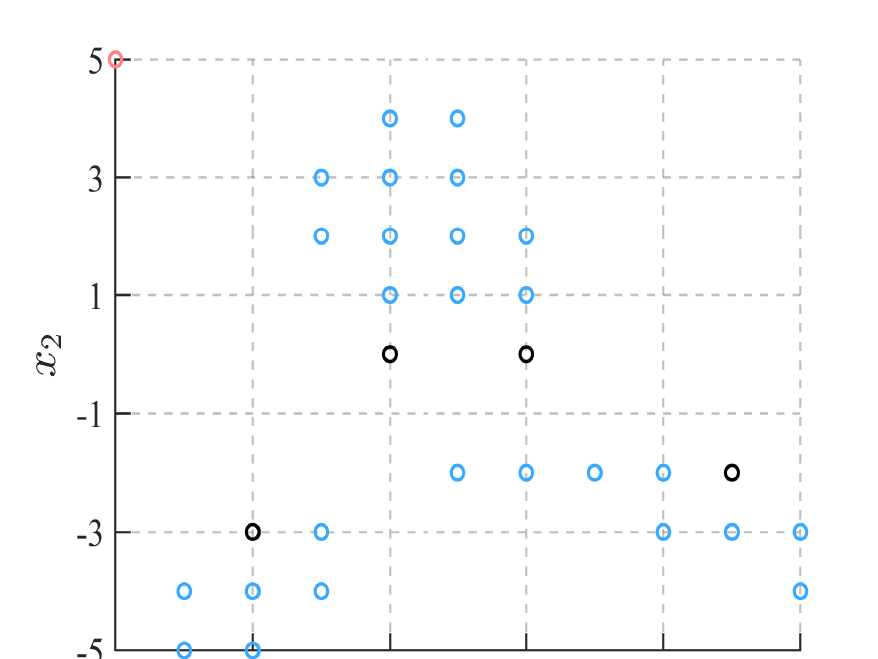}}
    \end{minipage}\ 
    \begin{minipage}[t]{0.49\linewidth}
        \centering
{\includegraphics[width=0.8\textwidth]{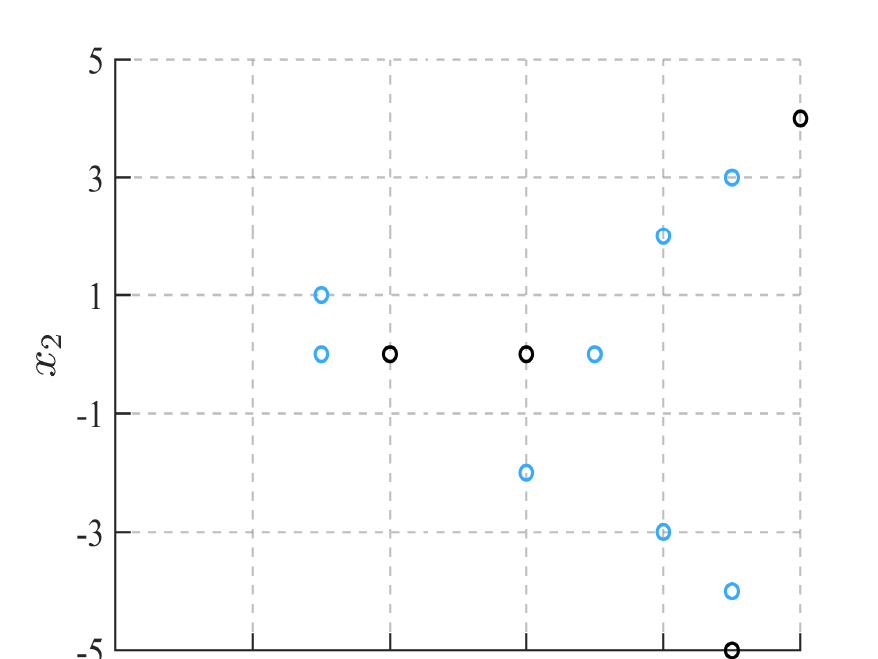}}
    \end{minipage}
    \caption{Illustration examples related to Conjecture \ref{thm5.1}. The bad point under  $\Lambda$-poisedness is marked as light blue, and the bad point under LOF as light red. The $y^3$ and $y^4$ in these two figures are both outlier traps. They are far away from $y^1$ and $y^2$ while having at least one light-blue point around them.}
        \label{fig:outlier-trap}
\end{figure}

During our experiments, an interesting phenomenon was observed. We find that bad points selected by $\Lambda$-poisedness tend to be near the fixed interpolation points. That is to say, when there are no other interpolation points nearby, there is always at least one bad point at the neighbor of one interpolation point. Fig.  \ref{fig:outlier-trap} shows this phenomenon. $y^3$ and $y^4$ in both these two figure are far away from $y^1$ and $y^2$. Here we mean $y^1$ and $y^2$ are not in the circle centered at $y^3$ or $y^4$, with radius.$\sqrt{2}$. While there is at least one bad point under $\Lambda$-poisedness in this region. We call such a point a trap, which can attract bad points.

Furthermore, to prove that this is no accident, we counted all the experimental points in the above experiment that satisfied the condition. Results are shown in Table \ref{tab:outlier-trap_stat}. Over 99\% of all experimental points far away from $y^1$ and $y^2$ is a trap. The only 0.9\% left may also be a trap because the resolution of grids in our experiment is only 1. This means that there may be bad points not covered but matching the condition as well.

\begin{table}[htbp]
	\centering
 	\caption{Times of $y^3$, $y^4$ being a outlier and a outlier trap. Here the concept of outlier is only restricted among $y^1$, $y^2$, $y^3$ and $y^4$.}
	\label{tab:outlier-trap_stat}
 \setlength{\tabcolsep}{20pt} 
	\begin{tabular}{llll}
    \toprule
		              & outliers & traps & rate\\
    \midrule
		$y^3$          & 2952 & 2925 & 99.1\% \\
	  $y^4$          & 2952 & 2925 & 99.1\% \\
		total          & 2952 & 2925 & 99.1\% \\
    \bottomrule
	\end{tabular}
\end{table}

It's understandable that when two interpolation points are nearby, the bad point can be either of two instead of points near them. Besides, This phenomenon also explains why bad points selected by LOF always differ from those by $\lambda$-poisedness.

\begin{conjecture}\label{thm5.1}
     For $\forall\ y^4\in\mathcal{Y}\backslash\{y^1,y^2,y^3\}$, if $\forall\ i\in\{1, 2, 3\},\lVert y^i-y^4\rVert>2$, then there is at least one $y^5\in\{y\lvert\lVert y-y^4\rVert\le2,y\in\mathcal{Y}\}$, which is the bad point under  $\Lambda$-poisedness.
\end{conjecture}

\subsection{{Directional Extension}}

As we look deeper, we find that the above conclusion is only apparent when the distribution of interpolated points is relatively decentralized. In this case, the distribution of interpolated points has no obvious ``directionality'', and bad points tend to cluster near the interpolated points. However, as the interpolated points come closer to each other, the range of bad points under  $\Lambda$-poisedness tends to extend significantly beyond the interpolated points when there is strong directionality, e.g., when three of the points are covariant.

\begin{figure}[htbp]
    \centering
    \begin{minipage}[t]{0.49\linewidth}
        \centering
        
        {\includegraphics[width=0.8\textwidth]{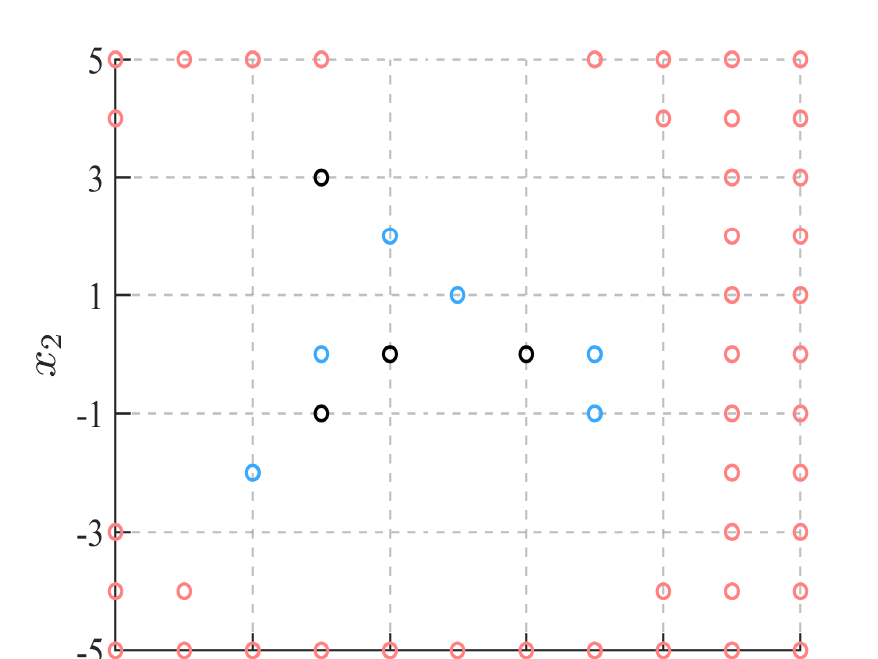}}
        
    \end{minipage}\ 
    \begin{minipage}[t]{0.49\linewidth}
        \centering
{\includegraphics[width=0.8\textwidth]{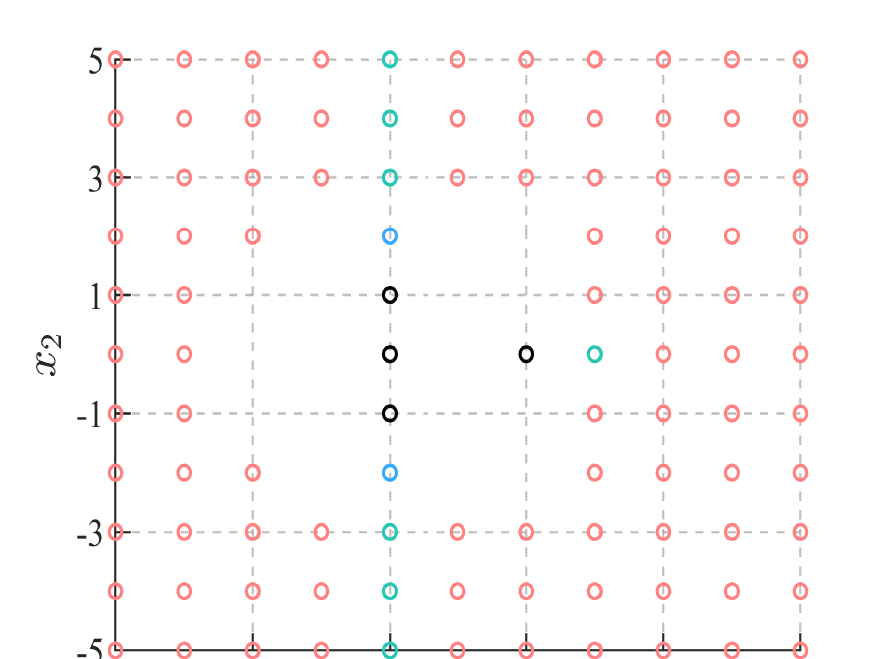}}
    \end{minipage}
    \caption{Illustration examples related to conjecture \ref{thm5.2}. Three-point congruence(left-hand side) versus another random condition (right-hand side). The bad point under  $\Lambda$-poisedness is marked blue, the bad point under LOF is marked  red, and points of both are marked green. }
        \label{fig:direction-ext}
\end{figure}

Fig.  \ref{fig:direction-ext} demonstrates the phenomena above. In these two figures, $y^1$, $y^2$, $y^3$ and $y^4$ are marked as black, the bad points under LOF as red, the bad points under $\Lambda$-poisedness as blue, and the points both bad under LOF and $\Lambda$-poisedness are highlighted as green. In the left panel of figure \ref{fig:direction-ext}, the location of all four points is \((-2, -1), (-1, 0), (1, 0), (-2, 3)\). Clearly, there is no strong ``directionality''. There are no three or more points on the same line, and none of the line segments connecting any two pairs of points is parallel to the other. As a result, all the blue points are near the black points, and there is no green point. However, in the right panel of Fig.  \ref{fig:direction-ext}, the location of the four points is $(-1, -1), (-1, 0), (-1, 1), (1, 0)$, where three points are in the same line. Then the bad points under  $\Lambda$-poisedness extend significantly beyond the interpolated points, just like the light green point in the left panel of Fig.  \ref{fig:direction-ext}. Here, we finally give the conjecture as follows.

\begin{conjecture}\label{thm5.2}
     For $\forall i,j,k\in\{1, 2, 3, 4\}, i\neq j\neq k$, if 
     \[
     \overrightarrow{y^i y^j}=\lambda \overrightarrow{y^j y^k},\]
     where $\lambda\in {\Re}, \lambda \neq 0$, there is at least one 
     \[y^5\in\{y\vert \lambda \overrightarrow{y^j y^k},\lambda\in {\Re}, \lambda \neq 0\},\]
     which is the  bad point under $\Lambda$-poisedness.
\end{conjecture}

\section{Conclusion}
\label{sec:conclusion}
In this study, we primarily discuss the concept of $\Lambda$-poisedness in the context of derivative-free optimization. $\Lambda$-poisedness describes the condition where, within a specified region, a set of interpolation points is well-suited for interpolating function values, ensuring that the condition number of the interpolation matrix is controlled by $\Lambda$. This property is essential in derivative-free optimization, as it allows us to rely on function values rather than directly obtaining gradient information. We compare $\Lambda$-poisedness with the Local Outlier Factor (LOF) used in data science. LOF measures the degree to which a data point is an outlier and is commonly employed in data mining and anomaly detection. We discuss the similarities and differences between $\Lambda$-poisedness and LOF. Through this study, we gain a deeper understanding of the role of $\Lambda$-poisedness in derivative-free optimization and highlight the significant impact of region selection on interpolation effectiveness. These insights are crucial for the design and improvement of optimization algorithms. At the same time, we proposed some conjectures.


\bibliography{sn-bibliography}


\begin{thebibliography}{16}
\ifx \bisbn   \undefined \def \bisbn  #1{ISBN #1}\fi
\ifx \binits  \undefined \def \binits#1{#1}\fi
\ifx \bauthor  \undefined \def \bauthor#1{#1}\fi
\ifx \batitle  \undefined \def \batitle#1{#1}\fi
\ifx \bjtitle  \undefined \def \bjtitle#1{#1}\fi
\ifx \bvolume  \undefined \def \bvolume#1{\textbf{#1}}\fi
\ifx \byear  \undefined \def \byear#1{#1}\fi
\ifx \bissue  \undefined \def \bissue#1{#1}\fi
\ifx \bfpage  \undefined \def \bfpage#1{#1}\fi
\ifx \blpage  \undefined \def \blpage #1{#1}\fi
\ifx \burl  \undefined \def \burl#1{\textsf{#1}}\fi
\ifx \doiurl  \undefined \def \doiurl#1{\url{https://doi.org/#1}}\fi
\ifx \betal  \undefined \def \betal{\textit{et al.}}\fi
\ifx \binstitute  \undefined \def \binstitute#1{#1}\fi
\ifx \binstitutionaled  \undefined \def \binstitutionaled#1{#1}\fi
\ifx \bctitle  \undefined \def \bctitle#1{#1}\fi
\ifx \beditor  \undefined \def \beditor#1{#1}\fi
\ifx \bpublisher  \undefined \def \bpublisher#1{#1}\fi
\ifx \bbtitle  \undefined \def \bbtitle#1{#1}\fi
\ifx \bedition  \undefined \def \bedition#1{#1}\fi
\ifx \bseriesno  \undefined \def \bseriesno#1{#1}\fi
\ifx \blocation  \undefined \def \blocation#1{#1}\fi
\ifx \bsertitle  \undefined \def \bsertitle#1{#1}\fi
\ifx \bsnm \undefined \def \bsnm#1{#1}\fi
\ifx \bsuffix \undefined \def \bsuffix#1{#1}\fi
\ifx \bparticle \undefined \def \bparticle#1{#1}\fi
\ifx \barticle \undefined \def \barticle#1{#1}\fi
\bibcommenthead
\ifx \bconfdate \undefined \def \bconfdate #1{#1}\fi
\ifx \botherref \undefined \def \botherref #1{#1}\fi
\ifx \url \undefined \def \url#1{\textsf{#1}}\fi
\ifx \bchapter \undefined \def \bchapter#1{#1}\fi
\ifx \bbook \undefined \def \bbook#1{#1}\fi
\ifx \bcomment \undefined \def \bcomment#1{#1}\fi
\ifx \oauthor \undefined \def \oauthor#1{#1}\fi
\ifx \citeauthoryear \undefined \def \citeauthoryear#1{#1}\fi
\ifx \endbibitem  \undefined \def \endbibitem {}\fi
\ifx \bconflocation  \undefined \def \bconflocation#1{#1}\fi
\ifx \arxivurl  \undefined \def \arxivurl#1{\textsf{#1}}\fi
\csname PreBibitemsHook\endcsname

\bibitem[\protect\citeauthoryear{Qian and Yu}{2021}]{Qian_Yu_2021}
\begin{barticle}
\bauthor{\bsnm{Qian}, \binits{H.}},
\bauthor{\bsnm{Yu}, \binits{Y.}}:
\batitle{Derivative-free reinforcement learning: a review}.
\bjtitle{Frontiers of Computer Science}
\bvolume{15},
\bfpage{156336}
(\byear{2021})
\end{barticle}
\endbibitem

\bibitem[\protect\citeauthoryear{Ghanbari and
  Scheinberg}{2017}]{Ghanbari_Scheinberg_2017}
\begin{botherref}
\oauthor{\bsnm{Ghanbari}, \binits{H.}},
\oauthor{\bsnm{Scheinberg}, \binits{K.}}:
Black-box optimization in machine learning with trust region based derivative
  free algorithm.
Technical Report 17T-005,
COR{@}L,
Bethlehem, PA, USA
(2017)
\end{botherref}
\endbibitem

\bibitem[\protect\citeauthoryear{Eldred et~al.}{2023}]{Eldred_Etal_2023}
\begin{barticle}
\bauthor{\bsnm{Eldred}, \binits{J.S.}},
\bauthor{\bsnm{Larson}, \binits{J.}},
\bauthor{\bsnm{Padidar}, \binits{M.}},
\bauthor{\bsnm{Stern}, \binits{E.}},
\bauthor{\bsnm{Wild}, \binits{S.M.}}:
\batitle{Derivative-free optimization of a rapid-cycling synchrotron}.
\bjtitle{Optimization and Engineering}
\bvolume{24},
\bfpage{1289}--\blpage{1319}
(\byear{2023})
\end{barticle}
\endbibitem

\bibitem[\protect\citeauthoryear{Conn
  et~al.}{2009}]{Conn_Scheinberg_Vicente_2009}
\begin{bbook}
\bauthor{\bsnm{Conn}, \binits{A.R.}},
\bauthor{\bsnm{Scheinberg}, \binits{K.}},
\bauthor{\bsnm{Vicente}, \binits{L.N.}}:
\bbtitle{Introduction to Derivative-Free Optimization}.
\bsertitle{MPS-SIAM Series on Optimization}.
\bpublisher{SIAM},
\blocation{Philadelphia, PA, USA}
(\byear{2009})
\end{bbook}
\endbibitem

\bibitem[\protect\citeauthoryear{Audet and Hare}{2017}]{Audet_Hare_2017}
\begin{bbook}
\bauthor{\bsnm{Audet}, \binits{C.}},
\bauthor{\bsnm{Hare}, \binits{W.}}:
\bbtitle{Derivative-Free and Blackbox Optimization}.
\bsertitle{Springer Series in Operations Research and Financial Engineering}.
\bpublisher{Springer},
\blocation{Cham, Switzerland}
(\byear{2017})
\end{bbook}
\endbibitem

\bibitem[\protect\citeauthoryear{Powell}{1998}]{Powell_1998}
\begin{barticle}
\bauthor{\bsnm{Powell}, \binits{M.J.D.}}:
\batitle{Direct search algorithms for optimization calculations}.
\bjtitle{Acta Numerica}
\bvolume{7},
\bfpage{287}--\blpage{336}
(\byear{1998})
\end{barticle}
\endbibitem

\bibitem[\protect\citeauthoryear{Cust{\'{o}}dio
  et~al.}{2017}]{Custodio_Scheinberg_Vicente_2017}
\begin{bchapter}
\bauthor{\bsnm{Cust{\'{o}}dio}, \binits{A.L.}},
\bauthor{\bsnm{Scheinberg}, \binits{K.}},
\bauthor{\bsnm{Vicente}, \binits{L.N.}}:
\bctitle{Methodologies and software for derivative-free optimization}.
In: \beditor{\bsnm{Terlaky}, \binits{T.}},
\beditor{\bsnm{Anjos}, \binits{M.F.}},
\beditor{\bsnm{Ahmed}, \binits{S.}} (eds.)
\bbtitle{Advances and Trends in Optimization with Engineering Applications}.
\bsertitle{MOS-SIAM Series on Optimization},
pp. \bfpage{495}--\blpage{506}.
\bpublisher{SIAM},
\blocation{Philadelphia, PA, USA}
(\byear{2017})
\end{bchapter}
\endbibitem

\bibitem[\protect\citeauthoryear{Larson
  et~al.}{2019}]{Larson_Menickelly_Wild_2019}
\begin{barticle}
\bauthor{\bsnm{Larson}, \binits{J.}},
\bauthor{\bsnm{Menickelly}, \binits{M.}},
\bauthor{\bsnm{Wild}, \binits{S.M.}}:
\batitle{Derivative-free optimization methods}.
\bjtitle{Acta Numerica}
\bvolume{28},
\bfpage{287}--\blpage{404}
(\byear{2019})
\end{barticle}
\endbibitem

\bibitem[\protect\citeauthoryear{Xie and x.~Yuan}{2023a}]{xie2023dfoto}
\begin{botherref}
\oauthor{\bsnm{Xie}, \binits{P.}},
\oauthor{\bsnm{Yuan}, \binits{Y.-x.}}:
Derivative-free optimization with transformed objective functions ({DFOTO}) and
  the algorithm based on the least {F}robenius norm updating quadratic model.
Journal of the Operations Research Society of China
(2023)
\end{botherref}
\endbibitem

\bibitem[\protect\citeauthoryear{Xie and x.~Yuan}{2023b}]{xie2023least}
\begin{botherref}
\oauthor{\bsnm{Xie}, \binits{P.}},
\oauthor{\bsnm{Yuan}, \binits{Y.-x.}}:
Least {$H^2$} norm updating quadratic interpolation model function for
  derivative-free trust-region algorithms
(2023)
\doiurl{10.48550/arXiv.2302.12017}
\end{botherref}
\endbibitem

\bibitem[\protect\citeauthoryear{Ragonneau and
  Zhang}{2024}]{ragonneau2024optimal}
\begin{botherref}
\oauthor{\bsnm{Ragonneau}, \binits{T.M.}},
\oauthor{\bsnm{Zhang}, \binits{Z.}}:
An optimal interpolation set for model-based derivative-free optimization
  methods.
Optimization Methods and Software,
1--13
(2024)
\end{botherref}
\endbibitem

\bibitem[\protect\citeauthoryear{Xie and x.~Yuan}{2023}]{xie2023linesearch}
\begin{barticle}
\bauthor{\bsnm{Xie}, \binits{P.}},
\bauthor{\bsnm{Yuan}, \binits{Y.-x.}}:
\batitle{A derivative-free optimization algorithm combining line-search and
  trust-region techniques}.
\bjtitle{Chinese Annals of Mathematics, Series B}
\bvolume{44}(\bissue{5}),
\bfpage{719}--\blpage{734}
(\byear{2023})
\end{barticle}
\endbibitem

\bibitem[\protect\citeauthoryear{Breunig et~al.}{2000}]{breunig_lof_2000}
\begin{bchapter}
\bauthor{\bsnm{Breunig}, \binits{M.M.}},
\bauthor{\bsnm{Kriegel}, \binits{H.-P.}},
\bauthor{\bsnm{Ng}, \binits{R.T.}},
\bauthor{\bsnm{Sander}, \binits{J.}}:
\bctitle{{LOF}: identifying density-based local outliers}.
In: \bbtitle{Proceedings of the 2000 {ACM} {SIGMOD} {International}
  {Conference} on {Management} of {Data}}.
\bsertitle{{SIGMOD} '00},
pp. \bfpage{93}--\blpage{104}.
\bpublisher{Association for Computing Machinery},
\blocation{New York, NY, USA}
(\byear{2000}).
\bcomment{event-place: Dallas, Texas, USA}
\end{bchapter}
\endbibitem

\bibitem[\protect\citeauthoryear{Tu et~al.}{2018}]{tu_hyperspectral_2018}
\begin{barticle}
\bauthor{\bsnm{Tu}, \binits{B.}},
\bauthor{\bsnm{Zhou}, \binits{C.}},
\bauthor{\bsnm{Kuang}, \binits{W.}},
\bauthor{\bsnm{Guo}, \binits{L.}},
\bauthor{\bsnm{Ou}, \binits{X.}}:
\batitle{Hyperspectral {Imagery} {Noisy} {Label} {Detection} by {Spectral}
  {Angle} {Local} {Outlier} {Factor}}.
\bjtitle{IEEE Geoscience and Remote Sensing Letters}
\bvolume{15}(\bissue{9}),
\bfpage{1417}--\blpage{1421}
(\byear{2018}).
\bcomment{Conference Name: IEEE Geoscience and Remote Sensing Letters}
\end{barticle}
\endbibitem

\bibitem[\protect\citeauthoryear{Peng et~al.}{2021}]{peng_electricity_2021}
\begin{barticle}
\bauthor{\bsnm{Peng}, \binits{Y.}},
\bauthor{\bsnm{Yang}, \binits{Y.}},
\bauthor{\bsnm{Xu}, \binits{Y.}},
\bauthor{\bsnm{Xue}, \binits{Y.}},
\bauthor{\bsnm{Song}, \binits{R.}},
\bauthor{\bsnm{Kang}, \binits{J.}},
\bauthor{\bsnm{Zhao}, \binits{H.}}:
\batitle{Electricity {Theft} {Detection} in {AMI} {Based} on {Clustering} and
  {Local} {Outlier} {Factor}}.
\bjtitle{IEEE Access}
\bvolume{9},
\bfpage{107250}--\blpage{107259}
(\byear{2021}).
\bcomment{Conference Name: IEEE Access}
\end{barticle}
\endbibitem

\bibitem[\protect\citeauthoryear{Wang and Lu}{2011}]{wang_efficient_2011}
\begin{barticle}
\bauthor{\bsnm{Wang}, \binits{W.}},
\bauthor{\bsnm{Lu}, \binits{P.}}:
\batitle{An {Efficient} {Switching} {Median} {Filter} {Based} on {Local}
  {Outlier} {Factor}}.
\bjtitle{IEEE Signal Processing Letters}
\bvolume{18}(\bissue{10}),
\bfpage{551}--\blpage{554}
(\byear{2011}).
\bcomment{Conference Name: IEEE Signal Processing Letters}
\end{barticle}
\endbibitem

\end{thebibliography}

\end{document}